\documentclass[a4paper,12pt]{article}
\usepackage{amsthm,amsmath,amssymb,latexsym,bm}

\newtheorem{theorem}{Theorem}[section]
\newtheorem{lemma}[theorem]{Lemma}
\theoremstyle{definition}

\theoremstyle{remark}
\newtheorem{remark}[theorem]{Remark}
\numberwithin{equation}{section}
\title{A class of higher order Painlev\'{e} systems arising from integrable hierarchies of type $A$}
\author{Takao Suzuki \thanks{Department of Mathematics, Kinki University, 3-4-1, Kowakae, Higashi-Osaka, Osaka 577-8502, Japan. E-mail: suzuki@math.kindai.ac.jp}}
\date{}

\begin{document}

\maketitle

% Abstract
\begin{abstract}
A relationship between Painlev\'{e} systems and infinite-dimensional integrable hierarchies is studied.
We derive a class of higher order Painlev\'{e} systems from Drinfeld-Sokolov (DS) hierarchies of type $A$ by similarity reductions.
This result allows us to understand some properties of Painlev\'{e} systems, Hamiltonian structures, Lax pairs and affine Weyl group symmetries.

\end{abstract}

% Section 1
\section{Introduction}

The connection between the second Painlev\'{e} equation and the KdV equation was clarified by Ablowitz and Segur \cite{AS}.
Since their result, a relationship between (higher order) Painlev\'{e} systems and infinite-dimensional integrable hierarchies has been studied.
In a recent work \cite{FS2}, a class of fourth order Painlev\'{e} systems was derived from the DS hierarchies of type $A$ by similarity reductions.
In this article, we give its development, namely, we derive a class of higher order Painlev\'{e} systems.

The DS hierarchies are extensions of the KdV hierarchy for the affine Lie algebras \cite{DS,GHM}.
They are characterized by the Heisenberg subalgebras of the affine Lie algebras.
And the isomorphism classes of the Heisenberg subalgebras are in one-to-one correspondence with the conjugacy classes of the finite Weyl group \cite{KP}.
Thus we can classify the DS hierarchies of type $A^{(1)}_n$ in terms of the partitions of the natural number $n+1$.
By means of this viewpoint, we list the known connections between Painlev\'{e} systems and integrable hierarchies of type $A$ in Table 1 and 2.
\begin{equation*}\begin{tabular}{|c|c|c|c|c|c|c|c|}
	\multicolumn{8}{c}{{\bf Table 1.} Painlev\'{e} equations and DS hierarchy}\\\hline
	Partition& $(2)$& $(1,1)$& $(3)$& $(2,1)$& $(1,1,1)$& $(4)$& $(2,2)$\\\hline
	Painlev\'{e} eq.& $P_{\rm{II}}$& $P_{\rm{IV}}$& $P_{\rm{IV}}$& $P_{\rm{V}}$& $P_{\rm{VI}}$& $P_{\rm{V}}$& $P_{\rm{VI}}$\\\hline
	Ref.& \cite{AS}& \cite{KK1}& \cite{Adl}& \cite{KIK}& \cite{KK2}& \cite{Adl}& \cite{FS2}\\\hline
\end{tabular}\end{equation*}
\begin{equation*}\begin{tabular}{|c|c|c|c|c|c|}
	\multicolumn{6}{c}{{\bf Table 2.} Higher order Painlev\'{e} systems and DS hierarchy}\\\hline
	Partition& $(3,1)$& $(4,1)$& $(2,2,1)$& $(3,3)$& $(n+1)$ for $n\geq4$\\\hline
	Painlev\'{e} sys.& $P_{(5)}$& $P_{(6)}$& $P_{(3,3)}$& $P_{(3,3)}$& $P_{(n+1)}$\\\hline
	Order of sys.& 4& 4& 4& 4& $n$ for $n$:even\\
	& & & & & $n-1$ for $n$:odd\\\hline
	Ref.& \cite{FS2}& \cite{FS2}& \cite{FS2}& \cite{FS2}& \cite{Adl,NY1}\\\hline
\end{tabular}\end{equation*}
Here the symbol $P_{(n+1)}$ stands for the higher order Painlev\'{e} system of type $A^{(1)}_{n}$ \cite{NY1}, or equivalently, the $(n+1)$-periodic Darboux chain \cite{Adl}.
The symbol $P_{(3,3)}$ stands for the fourth order Painlev\'{e} system with the coupled sixth Painlev\'{e} Hamiltonian \cite{FS2}; we describe its explicit formula below.

In this article, we consider a higher order generalization of the above facts.
The obtained results are listed in Table 3.
\begin{equation*}\begin{tabular}{|c|c|c|c|c|}
	\multicolumn{5}{c}{{\bf Table 3.} The result obtained in this article}\\\hline
	Partition& $(2n-1,1)$& $(2n,1)$& $(n,n,1)$& $(n+1,n+1)$\\\hline
	Painlev\'{e} sys.& $P_{(2n+1)}$& $P_{(2n+2)}$& $P_{(n+1,n+1)}$& $P_{(n+1,n+1)}$\\\hline
	Order of sys.& $2n$& $2n$& $2n$& $2n$\\\hline
	Ref.& \text{App.\ref{Sec:Other}}& \text{App.\ref{Sec:Other}}& \text{App.\ref{Sec:Other}}& \text{Sec.\ref{Sec:n+1n+1}}\\\hline
\end{tabular}\end{equation*}
The Painlev\'{e} system $P_{(n+1,n+1)}$ is a Hamiltonian system
\begin{equation*}
	\frac{dq_i}{dt} = \frac{\partial H}{\partial p_i},\quad
	\frac{dp_i}{dt} = -\frac{\partial H}{\partial q_i}\quad (i=1,\ldots,n),
\end{equation*}
with a coupled sixth Painlev\'{e} Hamiltonian
\begin{equation*}\begin{split}
	t(t-1)H &= \sum_{i=1}^{n}H_{\rm{VI}}\left[\sum_{j=0}^{n}\alpha_{2j+1}-\alpha_{2i-1}-\eta,\sum_{j=0}^{i-1}\alpha_{2j},\sum_{j=i}^{n}\alpha_{2j},\alpha_{2i-1}\eta;q_i,p_i\right]\\
	&\quad + \sum_{1\leq i<j\leq n}(q_i-1)(q_j-t)\{(q_ip_i+\alpha_{2i-1})p_j+p_i(p_jq_j+\alpha_{2j-1})\},
\end{split}\end{equation*}
where
\begin{equation*}\begin{split}
	H_{\rm{VI}}[\kappa_0,\kappa_1,\kappa_t,\kappa;q,p] &= q(q-1)(q-t)p^2 - \kappa_0(q-1)(q-t)p\\
	&\quad - \kappa_1q(q-t)p - (\kappa_t-1)q(q-1)p + \kappa q.
\end{split}\end{equation*}
Here the parameters $\alpha_0,\ldots,\alpha_{2n+1}$ satisfy the relation $\sum_{i=0}^{2n+1}\alpha_i=1$.
This system has a Lax pair associated with the loop algebra $\mathfrak{sl}_{2n+2}[z,z^{-1}]$ (or $\mathfrak{gl}_{2n+2}[z,z^{-1}]$) and admits the affine Weyl group symmetry of type $A^{(1)}_{2n+1}$; we discuss their details in Section \ref{Sec:CP6}.
Note that the system $P_{(2,2)}$ is equivalent to the sixth Painlev\'{e} equation.

\begin{remark}
The regular conjugacy classes of $W(A_n)$ correspond to the partitions $(p,\ldots,p)$ and $(p,\ldots,p,1)${\rm;} cf. {\rm\cite{DF,FHM}}.
Therefore any hierarchy in Table 1, 2 and 3 is associated with the regular conjugacy class of $W(A_n)$.
\end{remark}

\begin{remark}
The DS hierarchy for the partition $(n+1,n+1)$ is equivalent to the $(n+1,n+1)$-periodic reduction of the two-component KP hierarchy{\rm;} cf. {\rm\cite{BK,UT}}.
\end{remark}

\begin{remark}
The system $P_{(n+1,n+1)}$ is independently given by Tsuda as a similarity reduction of the UC hierarchy {\rm\cite{T1,T2}}.
Furthermore, $P_{(3,3)}$ appears in the classification of four-dimensional Painlev\'{e} type differential equations, which is given by Sakai {\rm\cite{Sak}}.
In both of them, the Painlev\'{e} system is given as the monodromy preserving deformation of a Fuchsian system, which is different from the Lax pair given in Section \ref{Sec:CP6} in terms of singularities and residue matrices.
The relationship between those two linear systems has been clarified with the aid of a Laplace transformation {\rm\cite{FS3}}.
\end{remark}

\begin{remark}
The system $P_{(n+1,n+1)}$ has a particular solution in terms of the generalized hypergeometric function ${}_{n+1}F_n$ {\rm\cite{Su,T3}}.
\end{remark}

\begin{remark}
The higher order Painlev\'{e} system of type $D_{2n+2}^{(1)}$, which is expressed as a Hamiltonian system of 2n-th order with a coupled sixth Painlev\'{e} Hamiltonian, was proposed by Sasano with the aid of algebraic geometry for initial value space {\rm\cite{Sas}}.
It is also derived from the Drinfeld-Sokolov hierarchy of type $D_{2n+2}^{(1)}$ by a similarity reduction {\rm\cite{FS1}}.
\end{remark}

This article is organized as follows.
In Section \ref{Sec:DS}, we first recall the affine Lie algebra of type $A^{(1)}_n$.
We next formulate the DS hierarchies of type $A^{(1)}_n$ and their similarity reductions.
In Section \ref{Sec:n+1n+1}, we consider the partition $(n+1,n+1)$ and derive the system $P_{(n+1,n+1)}$.
In Section \ref{Sec:CP6}, we discuss some properties of $P_{(n+1,n+1)}$, a Lax pair and a group of symmetries.

% Section 2
\section{DS hierarchy}\label{Sec:DS}

In this section, we first recall the affine Lie algebra of type $A^{(1)}_n$, following the notation in \cite{Kac,FS2}.
We next formulate the DS hierarchies of type $A^{(1)}_n$ and their similarity reductions.

% Section 2.1
\subsection{Affine Lie algebra}

The affine Lie algebra $\widehat{\mathfrak{g}}=\mathfrak{g}(A^{(1)}_n)$ is a Kac-Moody Lie algebra whose generalized Cartan matrix $A=\left[a_{i,j}\right]_{i,j=0}^{n}$ is defined by
\begin{equation*}\begin{array}{llll}
	a_{i,i}=2& (i=0,\ldots,n),\\[4pt]
	a_{i,i+1}=a_{n,0}=a_{i+1,i}=a_{0,n}=-1& (i=0,\ldots,n-1),\\[4pt]
	a_{i,j}=0& (\text{otherwise}).
\end{array}\end{equation*}
It is generated by the Chevalley generators $e_i,f_i,\alpha_i^{\vee}$ $(i=0,\ldots,n)$ and the scaling element $d$ with the fundamental relations
\begin{equation*}\begin{split}
	&[\alpha_i^{\vee},\alpha_j^{\vee}]=0,\quad
	[\alpha_i^{\vee},e_j]=a_{i,j}e_j,\quad
	[\alpha_i^{\vee},f_j]=-a_{i,j}f_j,\quad
	[e_i,f_j]=\delta_{i,j}\alpha_i^{\vee},\\
	&[d,\alpha_i^{\vee}]=0,\quad
	[d,e_i]=\delta_{i,0}e_0,\quad
	[d,f_i]=-\delta_{i,0}f_0,\\
	&(\mathrm{ad}e_i)^{1-a_{i,j}}(e_j)=0,\quad
	(\mathrm{ad}f_i)^{1-a_{i,j}}(f_j)=0\quad (i\neq j),
\end{split}\end{equation*}
for $i,j=0,\ldots,n$.
The canonical central element of $\widehat{\mathfrak{g}}$ is given by
\begin{equation*}
	K = \alpha_0^{\vee} + \alpha_1^{\vee} + \ldots + \alpha_n^{\vee}.
\end{equation*}
The normalized invariant form is given by the conditions
\begin{equation*}\begin{array}{lll}
	(\alpha_i^{\vee}|\alpha_j^{\vee}) = a_{i,j},& (e_i|f_j) = \delta_{i,j},&
	(\alpha_i^{\vee}|e_j) = (\alpha_i^{\vee}|f_j) = 0,\\[4pt]
	(d|d) = 0,& (d|\alpha_j^{\vee}) = \delta_{0,j},& (d|e_j) = (d|f_j) = 0,
\end{array}\end{equation*}
for $i,j=0,\ldots,n$.
We set
\begin{equation*}
	e_{i,j} = \mathrm{ad}e_i\mathrm{ad}e_{i+1}\ldots\mathrm{ad}e_{i+j-1}(e_{i+j}),\quad
	f_{i,j} = \mathrm{ad}f_{i+j}\mathrm{ad}f_{i+j-1}\ldots\mathrm{ad}f_{i+1}(f_i),
\end{equation*}
where $e_{i+n+1}=e_i$ and $f_{i+n+1}=f_i$.

The Cartan subalgebra of $\widehat{\mathfrak{g}}$ is defined by
\begin{equation*}
	\mathfrak{h} = \mathbb{C}\alpha_0^{\vee}\oplus\mathbb{C}\alpha_1^{\vee}\oplus\cdots\oplus\mathbb{C}\alpha_n^{\vee}\oplus\mathbb{C}d.
\end{equation*}
Let $\mathfrak{n}_{+}$ and $\mathfrak{n}_{-}$ be the subalgebras of $\widehat{\mathfrak{g}}$ generated by $e_i$ and $f_i$ $(i=0,\ldots,n)$, respectively.
Then the Borel subalgebra $\mathfrak{b}_{+}$ of $\widehat{\mathfrak{g}}$ is given by $\mathfrak{b}_{+}=\mathfrak{h}\oplus\mathfrak{n}_{+}$.
Note that we have the triangular decomposition
\begin{equation*}
	\widehat{\mathfrak{g}} = \mathfrak{n}_{-}\oplus\mathfrak{h}\oplus\mathfrak{n}_{+}
	= \mathfrak{n}_{-}\oplus\mathfrak{b}_{+}.
\end{equation*}

The isomorphism classes of the Heisenberg subalgebras of $\widehat{\mathfrak{g}}$ are in one-to-one correspondence with partitions of the natural number $n+1$.
Let $\mathbf{n}=(n_1,\ldots,n_k)$ be a partition of $n+1$.
Then the corresponding Heisenberg subalgebra is defined by
\begin{equation*}
	\mathfrak{s}_\mathbf{n} = \mathcal{P}_{n_1-1} \oplus \ldots \oplus \mathcal{P}_{n_k-1} \oplus \mathcal{H}_{k-1} \oplus \mathbb{C}K,
\end{equation*}
where $\mathcal{P}_n\oplus\mathbb{C}K$ and $\mathcal{H}_n\oplus\mathbb{C}K$ are isomorphic to the principal and homogeneous Heisenberg subalgebra of $\mathfrak{g}(A^{(1)}_n)$, respectively \cite{DF}.

The partition $\mathbf{n}$ determines a grading operator $\vartheta_{\mathbf{n}}\in\mathfrak{h}$, whose explicit formula is not given here (see Section 3 of \cite{FS2}).
The operator $\vartheta_{\mathbf{n}}$ defines a $\mathbb{Z}$-gradation of type $\mathbf{s}$ by
\begin{equation*}
	\widehat{\mathfrak{g}} = \bigoplus_{k\in\mathbb{Z}}\mathfrak{g}_k(\mathbf{s}),\quad
	\mathfrak{g}_k(\mathbf{s}) = \left\{x\in\widehat{\mathfrak{g}}\bigm|[\vartheta_{\mathbf{n}},x]=kx\right\},
\end{equation*}
where $\mathbf{s}=(s_0,\ldots,s_n)$ is a vector of non-negative integers given by
\begin{equation*}
	(\vartheta_{\mathbf{n}}|\alpha_i^{\vee}) = s_i\quad (i=0,\ldots,n).
\end{equation*}
Note that 
\begin{equation*}
	[\vartheta_{\mathbf{n}},e_i] = s_ie_i,\quad
	[\vartheta_{\mathbf{n}},f_i] = -s_if_i\quad (i=0,\ldots,n).
\end{equation*}
The Heisenberg subalgebra $\mathfrak{s}_\mathbf{n}$ admits the gradation defined by $\vartheta_{\mathbf{n}}$.

% Section 2.2
\subsection{DS hierarchy and similarity reduction}

The positive part of the Heisenberg subalgebra $\mathfrak{s}_\mathbf{n}$ has a graded basis $\{\Lambda_k\}_{k\in\mathbb{N}}$ satisfying
\begin{equation*}
	[\Lambda_k,\Lambda_l] = 0,\quad
	[\vartheta_{\mathbf{n}},\Lambda_k] = d_k\Lambda_k\quad (k,l\in\mathbb{N}),
\end{equation*}
where $d_k$ is a positive integer.
We assume that $d_k\leq d_{k+1}$ for any $k\in\mathbb{N}$.
In this subsection, we formulate the DS hierarchy associated with $\mathfrak{s}_\mathbf{n}$ by using those $\Lambda_k$.

Introducing time variables $t_k$ $(k\in\mathbb{N})$, we consider the Sato equation for an $\mathfrak{n}_{-}$-valued function $W=W(t_1,t_2,\ldots)$
\begin{equation}\label{Eq:Sato_DS}
	\partial_k - B_k = \exp(\mathrm{ad}W)(\partial_k-\Lambda_k)\quad (k\in\mathbb{N}),
\end{equation}
where $\partial_k=\partial/\partial t_k$ and $B_k$ stands for the $b_{+}$-component of $\exp(\mathrm{ad}W)(\Lambda_k)$.
The compatibility condition of \eqref{Eq:Sato_DS} gives the DS hierarchy
\begin{equation}\label{Eq:DS}
	[\partial_k-B_k,\partial_l-B_l] = 0\quad (k,l\in\mathbb{N}).
\end{equation}

We now require a similarity condition
\begin{equation}\label{Eq:Sato_DS_SR}
	\vartheta_{\mathbf{n}} - \rho - \sum_{k=1}^{\infty}d_kt_k\partial_k = \exp(\mathrm{ad}W)\left(\vartheta_{\mathbf{n}}-\rho-\sum_{k=1}^{\infty}d_kt_k\partial_k\right),
\end{equation}
with an element $\rho\in\mathfrak{h}$ satisfying
\begin{equation*}
	[\partial_k,\rho] = 0,\quad
	[\Lambda_k,\rho] = 0\quad (k\in\mathbb{N}).
\end{equation*}
Then the compatibility condition of \eqref{Eq:Sato_DS} and \eqref{Eq:Sato_DS_SR} gives
\begin{equation}\label{Eq:DS_SR}
	[\mathcal{M},\partial_k-B_k] = 0,\quad
	[\partial_k-B_k,\partial_l-B_l] = 0\quad (k,l\in\mathbb{N}),
\end{equation}
where
\begin{equation*}
	\mathcal{M} = \vartheta_{\mathbf{n}} - \rho - \sum_{k=1}^{\infty}d_kt_kB_k.
\end{equation*}
We call the system \eqref{Eq:DS_SR} a similarity reduction of the DS hierarchy.
Note that $\mathcal{M}$ is the $\mathfrak{b}_{+}$-component of $\exp(\mathrm{ad}W)(\vartheta_{\mathbf{n}}-\rho-\sum_{k=1}^{\infty}d_kt_k\Lambda_k)$.

In the following section, we always assume that $t_2=1$ and $t_k=0$ for any $k\geq3$.
Under this specialization, the system \eqref{Eq:DS_SR} is described as a system of ordinary differential equations
\begin{equation}\label{Eq:DS_SR_ODE}
	[\mathcal{M},\partial_1-B_1] = 0,\quad
	\mathcal{M} = \vartheta_{\mathbf{n}} - \rho - d_1t_1B_1 - d_2B_2,
\end{equation}
from which the Painlev\'{e} system is derived.

% Section 3
\section{Derivation of the system $P_{(n+1,n+1)}$}\label{Sec:n+1n+1}

In this section, we derive the system $P_{(n+1,n+1)}$ from the similarity reduction \eqref{Eq:DS_SR} for $\mathbf{n}=(n+1,n+1)$.
Here we set
\begin{equation*}\begin{split}
	&e_{i+2n+2} = e_i,\quad
	f_{i+2n+2} = f_i,\quad
	\alpha^{\vee}_{i+2n+2} = \alpha^{\vee}_i,\quad
	w_{i+2n+2} = w_i\\
	&\kappa_{i+2n+2} = \kappa_i,\quad
	\varphi_{i+2n+2} = \varphi_i,\quad
	u_{i+2n+2} = u_i,\quad
	v_{i+2n+2} = v_i.
\end{split}\end{equation*}

% Section 3.1
\subsection{Similarity reduction of the DS hierarchy}

At first, we give an explicit formula of the Heisenberg subalgebra $\mathfrak{s}_{(n+1,n+1)}$ of $\widehat{\mathfrak{g}}=\mathfrak{g}(A^{(1)}_{2n+1})$ following \cite{BK,FS2,KL}.
Let
\begin{equation*}\begin{split}
	&\Lambda_{2k-1} = \sum_{i=0}^{n}e_{2i+1,2k-1},\quad
	\Lambda_{2k} = \sum_{i=0}^{n}e_{2i+2,2k-1},\\
	&\bar{\Lambda}_{2k-1} = \sum_{i=0}^{n}f_{2i+1,2k-1},\quad
	\bar{\Lambda}_{2k} = \sum_{i=0}^{n}f_{2i+2,2k-1},
\end{split}\end{equation*}
for $k\in\mathbb{N}$.
Then $\mathfrak{s}_{(n+1,n+1)}$ is expressed as
\begin{equation*}
	\mathfrak{s}_{(n+1,n+1)} = \bigoplus_{k\in\mathbb{N}\setminus(2n+2)\mathbb{N}}\mathbb{C}\bar{\Lambda}_k \oplus \mathbb{C}K \oplus \bigoplus_{k\in\mathbb{N}\setminus(2n+2)\mathbb{N}}\mathbb{C}\Lambda_k.
\end{equation*}
The grading operator $\vartheta_{(n+1,n+1)}$ is given by
\begin{equation*}
	\vartheta_{(n+1,n+1)} = (n+1)d + \sum_{i=0}^{n}i(n-i+1)\alpha^{\vee}_{2i} + \sum_{i=0}^{n}\frac{(2i+1)n-2i^2}{2}\alpha^{\vee}_{2i+1}.
\end{equation*}
It implies a $\mathbb{Z}$-gradation of type $(1,0,\ldots,1,0)$, namely
\begin{equation*}
	(\vartheta_{(n+1,n+1)}|\alpha^{\vee}_{2i}) = 1,\quad
	(\vartheta_{(n+1,n+1)}|\alpha^{\vee}_{2i+1}) = 0,
\end{equation*}
for $i=0,\ldots,n$.
Note that
\begin{equation*}
	[\vartheta_{(n+1,n+1)},\Lambda_{2k-1}] = k\Lambda_{2k-1},\quad
	[\vartheta_{(n+1,n+1)},\Lambda_{2k}] = k\Lambda_{2k}\quad (k\in\mathbb{N}).
\end{equation*}

The similarity reduction \eqref{Eq:DS_SR_ODE} for $\mathbf{n}=(n+1,n+1)$ is described as
\begin{equation}\label{Eq:SR_n+1n+1_1}
	[\mathcal{M},\partial_1-B_1] = 0,\quad
	\mathcal{M} = \vartheta_{(n+1,n+1)} - \rho - t_1B_1 - B_2,
\end{equation}
where
\begin{equation*}
	\rho = \rho_1\sum_{i=0}^{n}\alpha^{\vee}_{2i+1},\quad
	[\partial_1,\rho] = 0.
\end{equation*}
Let us to denote the $\mathfrak{b}_{+}$-valued functions $\mathcal{M}$ and $B_1$ by
\begin{equation*}
	\mathcal{M} = \kappa - \sum_{i=0}^{2n+1}\varphi_ie_i - t_1\Lambda_1 - \Lambda_2,\quad
	B_1 = u + \sum_{i=0}^{2n+1}v_ie_i + \Lambda_1,
\end{equation*}
where
\begin{equation*}
	\kappa = \vartheta_{(n+1,n+1)} - \sum_{i=0}^{2n+1}\kappa_i\alpha^{\vee}_i,\quad
	u = \sum_{i=0}^{2n+1}u_i\alpha^{\vee}_i.
\end{equation*}
Then the system \eqref{Eq:SR_n+1n+1_1} is rewritten into
\begin{equation}\label{Eq:SR_n+1n+1_2}
	\partial_1(\kappa_i) = 0,\quad
	\partial_1(\varphi_i) = (u|\alpha^{\vee}_i)\varphi_i + v_i(\kappa|\alpha^{\vee}_i),
\end{equation}
for $i=0,\ldots,2n+1$ and
\begin{equation}\begin{split}\label{Eq:SR_n+1n+1_3}
	&(u|\alpha^{\vee}_{2i}+\alpha^{\vee}_{2i+1}) - v_{2i+1}\varphi_{2i} + v_{2i}\varphi_{2i+1} = 0,\\
	&t_1(u|\alpha^{\vee}_{2i+1}+\alpha^{\vee}_{2i+2}) - v_{2i+2}\varphi_{2i+1} + v_{2i+1}\varphi_{2i+2} + (\kappa|\alpha^{\vee}_{2i+1}+\alpha^{\vee}_{2i+2}) = 1,\\
	&t_1v_{2i} - v_{2i+2} - \varphi_{2i} = 0,\quad
	v_{2i+1} - t_1v_{2i+3} + \varphi_{2i+3} = 0,
\end{split}\end{equation}
for $i=0,\ldots,n$.
In the next subsection, we express the system \eqref{Eq:SR_n+1n+1_2} with \eqref{Eq:SR_n+1n+1_3} as a Hamiltonian system.

% Section 3.2
\subsection{Hamiltonian system}

Let us to denote the $\mathfrak{n}_{-}$-valued function $W$ by
\begin{equation*}
	W = -\sum_{i=0}^{2n+1}w_if_i - \sum_{k=1}^{\infty}\sum_{i=0}^{2n+1}w_{i,k}f_{i,k}.
\end{equation*}
In this subsection, we express the system \eqref{Eq:SR_n+1n+1_2} with \eqref{Eq:SR_n+1n+1_3} as the Hamiltonian system in terms of the dependent variables $w_{2i+1},\varphi_{2i+1}$ $(i=0,\ldots,n)$.
Note that those variables are taken from the $\mathfrak{g}_0(1,0,\ldots,1,0)$-components of $W$ and $\mathcal{M}$.

The function $\mathcal{M}$ is defined as the $\mathfrak{b}_{+}$-components of $\exp(\mathrm{ad}W)(\vartheta_{(n+1,n+1)}-\rho-t_1\Lambda_1-\Lambda_2)$, from which we obtain
\begin{equation}\begin{split}\label{Eq:SR_n+1n+1_4}
	&\kappa_{2i} = -\frac{t_1}{2}w_{2i-1}w_{2i} + t_1w_{2i-1,1} + \frac{1}{2}w_{2i}w_{2i+1} + w_{2i,1},\\
	&\kappa_{2i+1} =  -\frac{1}{2}w_{2i}w_{2i+1} + w_{2i,1} + \frac{t_1}{2}w_{2i+1}w_{2i+2} + t_1w_{2i+1,1} + \rho_1,\\
	&\varphi_{2i} = -t_1w_{2i-1} + w_{2i+1},\quad
	\varphi_{2i+1} = -w_{2i} + t_1w_{2i+2},
\end{split}\end{equation}
Note that
\begin{equation*}
	w_{2i+1} = -\sum_{j=0}^{n}\frac{t_1^j}{t_1^{n+1}-1}\varphi_{2i-2j}\quad (i=0,\ldots,n).
\end{equation*}
Similarly, the function $B_1$ is defined as the $\mathfrak{b}_{+}$-components of $\exp(\mathrm{ad}W)(\Lambda_1)$, from which we obtain
\begin{equation}\begin{split}\label{Eq:SR_n+1n+1_5}
	&u_{2i} = -\frac{1}{2}w_{2i-1}w_{2i} + w_{2i-1,1},\quad
	u_{2i+1} = \frac{1}{2}w_{2i+1}w_{2i+2} + w_{2i+1,1},\\
	&v_{2i} = -w_{2i-1},\quad
	v_{2i+1} = w_{2i+2}.
\end{split}\end{equation}
for $i=0,\ldots,n$.
Combining the equations \eqref{Eq:SR_n+1n+1_4} and \eqref{Eq:SR_n+1n+1_5}, we have
\begin{lemma}\label{Lem:SR_n+1n+1_Ham_Var}
The $\mathfrak{b}_{+}$-valued functions $\mathcal{M}$ and $B_1$ can be expressed in terms of the dependent variables $w_{2i+1},\varphi_{2i+1}$ $(i=0,\ldots,n)$ as
\begin{equation*}\begin{split}
	&\varphi_{2i} = -t_1w_{2i-1} + w_{2i+1},\\
	&u_{2i} - u_{2i+1} = -\sum_{j=0}^{n}\frac{t_1^{j-1}}{t_1^{n+1}-1}w_{2i+1}\varphi_{2i+2j+1} + \frac{1}{t_1}(\rho_1+\kappa_{2i}-\kappa_{2i+1}),\\
	&u_{2i+1} - u_{2i+2} = \sum_{j=0}^{n}\frac{t_1^j}{t_1^{n+1}-1}w_{2i+1}\varphi_{2i+2j+3},\\
	&v_{2i} = -w_{2i-1},\quad
	v_{2i+1} = \sum_{j=0}^{n}\frac{t_1^j}{t_1^{n+1}-1}\varphi_{2i+2j+3},
\end{split}\end{equation*}
for $i=1,\ldots,n+1$.
Furthermore, those variables satisfy
\begin{equation*}
	\sum_{i=0}^{n}w_{2i+1}\varphi_{2i+1} = -\sum_{i=0}^{n}(\rho_1+\kappa_{2i}-\kappa_{2i+1}).
\end{equation*}
\end{lemma}

Following \cite{NY2}, we define the Poisson structure for the function $\mathcal{M}$ by
\begin{equation*}
	\{\varphi_{2i},\varphi_{2i+1}\} = -(n+1),\quad
	\{\varphi_{2i+1},\varphi_{2i+2}\} = -(n+1)t_1\quad (i=0,\ldots,n).
\end{equation*}
Then we arrive at
\begin{theorem}
In terms of the variables $w_{2i+1},\varphi_{2i+1}$ $(i=0,\ldots,n)$ with the Poisson structure
\begin{equation*}
	\{\varphi_{2i+1},w_{2j+1}\} = (n+1)\delta_{i,j}\quad (i,j=0,\ldots,n),
\end{equation*}
the similarity reduction \eqref{Eq:SR_n+1n+1_1} is expressed as the Hamiltonian system
\begin{equation}\label{Eq:FS}
	\partial_1(w_{2i+1}) = \{H,w_{2i+1}\},\quad
	\partial_1(\varphi_{2i+1}) = \{H,\varphi_{2i+1}\}\quad (i=0,\ldots,n),
\end{equation}
with the Hamiltonian
\begin{equation}\begin{split}\label{Eq:FS_Ham}
	H &= \sum_{i=0}^{n}\frac{n}{2(n+1)^2t_1}(w_{2i+1}\varphi_{2i+1}+2\kappa_{2i}-2\kappa_{2i+1})w_{2i+1}\varphi_{2i+1}\\
	&\quad - \sum_{i=0}^{n}\sum_{j=1}^{n}\frac{1}{(n+1)^2t_1}(w_{2i+2j+1}\varphi_{2i+2j+1}+\kappa_{2i+2j}-\kappa_{2i+2j+1})w_{2i+1}\varphi_{2i+1}\\
	&\quad - \sum_{i=0}^{n}\sum_{j=0}^{n}\frac{t_1^j}{(n+1)(t_1^{n+1}-1)}\{w_{2i+1}\varphi_{2i+1}+(\kappa|\alpha^{\vee}_{2i+1})\}w_{2i+1}\varphi_{2i+2j+3},
\end{split}\end{equation}
and the relation
\begin{equation}\label{Eq:FS_var_rel}
	\sum_{i=0}^{n}w_{2i+1}\varphi_{2i+1} = -\sum_{i=0}^{n}(\rho_1+\kappa_{2i}-\kappa_{2i+1}).
\end{equation}
\end{theorem}

The system \eqref{Eq:FS} with \eqref{Eq:FS_Ham} and \eqref{Eq:FS_var_rel} can be rewritten into the Hamiltonian system in terms of the canonical coordinates.
The equation \eqref{Eq:FS_var_rel} implies
\begin{equation*}\begin{split}
	\sum_{i=0}^{n}d\varphi_{2i+1}\wedge dw_{2i+1} &= \sum_{i=0}^{n-1}d\varphi_{2i+1}\wedge dw_{2i+1} - \sum_{i=0}^{n-1}d\frac{w_{2i+1}\varphi_{2i+1}}{w_{2n+1}}\wedge dw_{2n+1}\\
	&= \sum_{i=0}^{n-1}d(w_{2n+1}\varphi_{2i+1})\wedge d\frac{w_{2i+1}}{w_{2n+1}}.
\end{split}\end{equation*}
Therefore we can take
\begin{equation}\label{Eq:FS_to_CP6}
	q_i = \frac{w_{2i-1}}{t_1^iw_{2n+1}},\quad
	p_i = \frac{t_1^iw_{2n+1}\varphi_{2i-1}}{n+1}\quad (i=1,\ldots,n),
\end{equation}
as canonical coordinates of a $2n$-dimensional system with a Poisson structure
\begin{equation*}
	\{p_i,q_j\} = \delta_{i,j}\quad (i,j=1,\ldots,n).
\end{equation*}
We denote the parameters by
\begin{equation*}
	\alpha_i = \frac{(\kappa|\alpha^{\vee}_i)}{n+1}\quad (i=0,\ldots,2n+1),\quad
	\eta = \sum_{j=0}^{n}\frac{\rho_1+\kappa_{2j}-\kappa_{2j+1}}{n+1}.
\end{equation*}
Via a transformation of the independent variable $t=t_1^{-(n+1)}$, we obtain
\begin{theorem}
The variables $q_i,p_i$ $(i=1,\ldots,n)$ defined by \eqref{Eq:FS_to_CP6} satisfy the Painlev\'{e} system $P_{(n+1,n+1)}$.
Then the variable $w_{2n+1}$ satisfies
\begin{equation*}\begin{split}
	t(t-1)\frac{d}{dt}\log w_{2n+1} &= -\sum_{i=1}^{n}\left\{(q_i-1)(q_i-t)p_i+\alpha_{2i-1}q_i\right\} - \alpha_{2n+1}\\
	&\quad + \frac{nt+n+2}{n+1}\eta + \sum_{i=0}^{n}\frac{n-2i}{2n+2}(\alpha_{2i-1}+\alpha_{2i})(t-1).
\end{split}\end{equation*}
\end{theorem}

We remark that the parameter $\eta$ satisfies
\begin{equation*}
	\{\eta,q_i\} = \{\eta,p_i\} = 0\quad (i=1,\ldots,n),\quad
	\{\eta,w_{2n+1}\} = w_{2n+1}.
\end{equation*}
Thus the Poisson algebra generated by $w_{2i+1},\varphi_{2i+1}$ $(i=0,\ldots,n)$ is equivalent to one generated by $q_i,p_i$ $(i=1,\ldots,n)$, $w_{2n+1}$ and $\eta$.

% Section 4
\section{Properties of the system $P_{(n+1,n+1)}$}\label{Sec:CP6}

In this section, we discuss some properties of the system $P_{(n+1,n+1)}$, a Lax pair and a group of symmetries.
In the following, we use a notation $q_{n+1}=t$ for a convenience.

% Section 4.1
\subsection{Lax pair}

In this subsection, we give a Lax pair of $P_{(n+1,n+1)}$ in a framework of the loop algebra $\mathfrak{gl}_{2n+2}[z,z^{-1}]$.

We consider a system of linear differential equations
\begin{equation}\label{Eq:SR_n+1n+1_Lax}
	z\frac{d\psi}{dz} = M\psi,\quad
	t\frac{d\psi}{dt} = B\psi.
\end{equation}
The matrix $M$ is given by
\begin{equation*}
	M = \begin{bmatrix}\varepsilon_1&\varphi_1&1\\&\varepsilon_2&\varphi_2&1\\&&&\ddots\\&&&&\varphi_{2n-1}&1\\&&&&\varepsilon_{2n}&\varphi_{2n}&1\\t^{-1}z&&&&&\varepsilon_{2n+1}&\varphi_{2n+1}\\\varphi_0z&z&&&&&\varepsilon_{2n+2}\end{bmatrix},
\end{equation*}
where
\begin{equation*}
	\varphi_{2i-1} = p_i,\quad
	\varphi_{2n+1} = -\frac{1}{t}\left(\sum_{j=1}^{n}q_jp_j+\eta\right),\quad
	\varphi_{2i} = q_{i+1} - q_i,\quad
	\varphi_0 = q_1 - 1,
\end{equation*}
for $i=1,\ldots,n$.
The matrix $B$ is given by
\begin{equation*}
	B = \begin{bmatrix}u_1&v_1&-1\\&u_2&v_2&0\\&&u_3&v_3&-1\\&&&u_4&v_4&0\\&&&&&\ddots\\&&&&&&v_{2n-1}&-1\\&&&&&&u_{2n}&v_{2n}&0\\-t^{-1}z&&&&&&&u_{2n+1}&v_{2n+1}\\v_0z&0&&&&&&&u_{2n+2}\end{bmatrix}.
\end{equation*}
where
\begin{equation*}\begin{split}
	&u_1 = tv_{2n+1}q_1 - \varepsilon_1 + \varepsilon_{2n+2} - 1,\quad
	u_{2i+1} = v_{2i-1}q_{i+1} - \varepsilon_{2i+1} + \varepsilon_{2n+2} - 1,\\
	&u_{2i} = -v_{2i-1}q_i + \frac{\sum_{i=1}^{n}\left\{(q_i-1)(q_i-t)p_i+\alpha_{2i-1}q_i\right\}+\alpha_{2n+1}t-(t+1)\eta}{t-1},\\
	&u_{2n+2} = -tv_{2n+1} + \frac{\sum_{i=1}^{n}\left\{(q_i-1)(q_i-t)p_i+\alpha_{2i-1}q_i\right\}+\alpha_{2n+1}t-(t+1)\eta}{t-1},\\
	&v_{2i-1} = -\frac{\sum_{j=1}^{i}(q_j-1)p_j+\sum_{j=i+1}^{n}(q_j-t)p_j+\eta}{t-1},\\
	&v_{2n+1} = -\frac{\sum_{j=1}^{n}(q_j-t)p_j+\eta}{t(t-1)},\quad
	v_{2i} = q_i,\quad
	v_0 = 1.
\end{split}\end{equation*}
for $i=1,\ldots,n$.
Then we have
\begin{theorem}
The compatibility condition of the system \eqref{Eq:SR_n+1n+1_Lax} gives the Painlev\'{e} system $P_{(n+1,n+1)}$ with the parameters
\begin{equation*}
	\alpha_0 = \varepsilon_1 - \varepsilon_{2n+2} + 1,\quad
	\alpha_i = -\varepsilon_i + \varepsilon_{i+1}\quad (i=1,\ldots 2n+1).
\end{equation*}
\end{theorem}

Such a Lax pair is given as follows.
Under the Sato equation \eqref{Eq:Sato_DS} and the similarity condition \eqref{Eq:Sato_DS_SR}, we consider a wave function
\begin{equation*}
	\Psi = \exp(W)z^{\rho}\exp\left(\sum_{k=1}^{\infty}t_k\Lambda_k\right).
\end{equation*}
Then we obtain a system of linear differential equations
\begin{equation}\label{Eq:DS_SR_Lax}
	\mathcal{M}\Psi = 0,\quad
	(\partial_k-B_k)\Psi\quad (k\in\mathbb{N}).
\end{equation}
Under the specialization $t_1=t^{-1/(n+1)}$, $t_2=1$ and $t_k=0$ $(k\geq3)$, the system \eqref{Eq:DS_SR_Lax} for $\mathbf{n}=(n+1,n+1)$ is transformed into the one \eqref{Eq:SR_n+1n+1_Lax} via a certain gauge transformation.

% Section 4.2
\subsection{Affine Weyl group symmetry}

The system $P_{(n+1,n+1)}$ admits an extended affine Weyl group symmetry of type $A^{(1)}_{2n+1}$.
In this subsection, we describe its action on the dependent variables and parameters.

At first, we define an extended affine Weyl group $\widetilde{W}(A^{(1)}_{2n+1})$.
It is generated by the transformations $r_0,\ldots,r_{2n-1}$ and $\pi$ with the fundamental relations
\begin{equation*}\begin{split}
	&r_i^2=1,\quad
	(r_ir_j)^{2-a_{i,j}}=1\quad (i,j=0,\ldots,2n+1; i\neq j),\\
	&\pi^{2n+2} = 1,\quad
	\pi r_i = r_{i+1}\pi,\quad
	\pi r_{2n+1} = r_0\pi\quad (i=0,\ldots,2n),
\end{split}\end{equation*}
where
\begin{equation*}\begin{array}{llll}
	a_{i,i}=2& (i=0,\ldots,2n+1),\\[4pt]
	a_{i,i+1}=a_{2n+1,0}=a_{i+1,i}=a_{0,2n+1}=-1& (i=0,\ldots,2n),\\[4pt]
	a_{i,j}=0& (\text{otherwise}).
\end{array}\end{equation*}

Let $r_0,\ldots,r_{2n+1}$ be birational canonical transformations defined by
\begin{equation*}\begin{split}
	&r_0(q_j) = q_j,\quad
	r_0(p_j) = p_j + \frac{\alpha_0}{q_1-1}\{q_1-1,p_j\},\\
	&r_{2i}(q_j) = q_j,\quad
	r_{2i}(p_j) = p_j + \frac{\alpha_{2i}}{q_i-q_{i+1}}\{q_i-q_{i+1},p_j\}\quad (i=1,\ldots,n),\\
	&r_{2i-1}(q_j) = q_j + \frac{\alpha_{2i-1}}{p_i}\{p_i,q_j\},\quad
	r_{2i-1}(p_j) = p_j\quad (i=1,\ldots,n),\\
	&r_{2n+1}(q_j) = q_j + \frac{\alpha_{2n+1}q_j}{\sum_{j=1}^{n}q_jp_j+\eta-\alpha_{2n+1}},\quad
	r_{2n+1}(p_j) = p_j - \frac{\alpha_{2n+1}p_j}{\sum_{j=1}^{n}q_jp_j+\eta},
\end{split}\end{equation*}
for $j=1,\ldots,n$ and
\begin{equation*}
	r_i(\alpha_j) = \alpha_j - a_{i,j}\alpha_i,\quad
	r_i(\eta) = \eta + (-1)^i\alpha_i\quad (i,j=0,\ldots,2n+1).
\end{equation*}
Also let $\pi$ be birational canonical transformations defined by
\begin{equation*}\begin{split}
	\pi(q_i) &= \frac{\sum_{j=1}^{i}(q_j-1)p_j+\sum_{j=i+1}^{n}(q_j-t)p_j+\eta}{\sum_{j=1}^{n}(q_j-t)p_j+\eta},\\
	\pi(p_i) &= \frac{(q_i-q_{i+1})\{\sum_{j=1}^{n}(q_j-1)p_j+\eta\}}{t-1}\quad (i=1,\ldots,n),
\end{split}\end{equation*}
and
\begin{equation*}
	\pi(t) = \frac{1}{t},\quad
	\pi(\alpha_i) = \alpha_{i+1}\quad (i=0,\ldots,2n),\quad
	\pi(\alpha_{2n+1}) = \alpha_0,\quad
	\pi(\eta) = -\eta.
\end{equation*}
Then we have
\begin{theorem}
The Painlev\'{e} system $P_{(n+1,n+1)}$ is invariant under actions of the transformations $r_0,\ldots,r_{2n-1}$ and $\pi$.
Furthermore, the group of symmetries $\langle r_0,\ldots,r_{2n-1},\pi\rangle$ is isomorphic to the extended affine Weyl group $\widetilde{W}(A^{(1)}_{2n+1})$.
\end{theorem}

Note that the transformations $r_0,\ldots,r_{2n+1}$ arise from gauge transformations for the system \eqref{Eq:DS_SR_Lax}
\[
	r_i(\Psi) = \exp\left(-\frac{\alpha_i}{\varphi_i}f_i\right)\Psi\quad (i=0,\ldots,2n+1).
\]

% Appendix
\appendix

% Section A
\section{Investigation for other partitions}\label{Sec:Other}

In the previous section, we derive the system $P_{(n+1,n+1)}$ from the DS hierarchy for the partition $(n+1,n+1)$.
In the same manner, we can derive the systems $P_{(n+1,n+1)}$, $P_{(2n+2)}$ and $P_{(2n+1)}$ from the hierarchies for the partitions $(n,n,1)$, $(2n,1)$ and $(2n-1,1)$ respectively.

We recall that the Painlev\'{e} systems $P_{(2n+2)}$ and $P_{(2n+1)}$ are Hamiltonian systems with the coupled Hamiltonians $H_{(2n+2)}$ and $H_{(2n+1)}$ respectively.
They are given by
\begin{equation*}\begin{split}
	tH_{(2n+2)} &= \sum_{i=1}^{n}H_{\rm{V}}\left[\alpha_{2i},\sum_{j=1}^{i}\alpha_{2j-1},\sum_{j=1}^{n+1}\alpha_{2j-1};q_i,p_i\right] + \sum_{1\leq i<j\leq n}2q_ip_i(q_j-1)p_j,\\
	H_{(2n+1)} &= \sum_{i=1}^{n}H_{\rm{IV}}\left[\alpha_{2i},\sum_{j=1}^{i}\alpha_{2j-1};q_i,p_i\right] + \sum_{1\leq i<j\leq n}2q_ip_ip_j,
\end{split}\end{equation*}
where
\begin{equation*}\begin{split}
	H_{\rm{V}}[a,b,c;q,p] &= q(q-1)p(p+t) + atq + bp - cqp,\\
	H_{\rm{IV}}[a,b;q,p] &= qp(p-q-t) - aq - bp.
\end{split}\end{equation*}

% Section A.1
\subsection{Partition $(n,n,1)$}

The Heisenberg subalgebra $\mathfrak{s}_{(n,n,1)}$ of $\mathfrak{g}(A_{2n}^{(1)})$ is generated by elements $\{\Lambda_k,\bar{\Lambda}_k\}_{k\in\mathbb{N}}$ with
\begin{equation*}
	\Lambda_1 = \sum_{i=1}^{n-1}e_{2i-1,1} + e_{2n-1,2},\quad
	\Lambda_2 = \sum_{i=1}^{n-1}e_{2i,1} + e_{2n,2}.
\end{equation*}
The grading operator $\vartheta_{(n,n,1)}$ is given by
\begin{equation*}
	\vartheta_{(n,n,1)} = nd + \sum_{i=1}^{n}(n-i+1)\left(i-\frac{n}{2n+1}\right)\alpha^{\vee}_{2i-1} + \sum_{i=1}^{n}i\left(\frac{2n}{2n+1}-i\right)\alpha^{\vee}_{2i}.
\end{equation*}
It implies a $\mathbb{Z}$-gradation of type $(0,1,0,\ldots,1,0)$.

The similarity reduction \eqref{Eq:DS_SR_ODE} for $\mathbf{n}=(n,n,1)$ is given by
\begin{equation}\label{Eq:SR_nn1}
	[\mathcal{M},\partial_1-B_1] = 0,
\end{equation}
with
\begin{equation*}
	\mathcal{M} = \kappa - \sum_{i=0}^{2n}\varphi_ie_i - \varphi_{0,1}e_{0,1} - \varphi_{2n-1,1}e_{2n-1,1} - \varphi_{2n,1}e_{2n,1} - t_1\Lambda_1 - \Lambda_2,
\end{equation*}
where
\begin{equation*}
	\kappa = \vartheta_{(n,n,1)} - \sum_{i=0}^{2n}\kappa_i\alpha^{\vee}_i.
\end{equation*}
Note that $\mathcal{M}$ is defined as a $\mathfrak{b}_{+}$-component of $\exp(\mathrm{ad}W)(\vartheta_{(n,n,1)}-\rho-t_1\Lambda_1-\Lambda_2)$, where
\begin{equation*}
	\rho = \rho_1\sum_{i=1}^{n}2i(\alpha^{\vee}_{2i-1}+\alpha^{\vee}_{2i}) + \rho_2\sum_{i=1}^{n}2(n-i+1)(\alpha^{\vee}_{2i-2}+\alpha^{\vee}_{2i-1}).
\end{equation*}

\begin{theorem}
Under the system \eqref{Eq:SR_nn1}, we set
\begin{equation*}
	t = t_1^{-n},\quad
	q_i = \frac{w_{2i-1}\varphi_{0,1}}{t_1^{i-1}(w_1\varphi_{0,1}-\varphi_0)},\quad
	p_i = \frac{t_1^{i-1}(w_1\varphi_{0,1}-\varphi_0)\varphi_{2i-1}}{n\varphi_{0,1}},
\end{equation*}
where
\[
	w_{2i-1} = -\sum_{j=1}^{i-1}\frac{t_1^{i-j-1}}{t_1^n-1}\varphi_{2j} - \sum_{j=i}^{n-1}\frac{t_1^{n+i-j-1}}{t_1^n-1}\varphi_{2j} - \frac{t_1^{i-1}}{t_1^n-1}\varphi_{2n,1},
\]
for $i=1,\ldots,n$ and
\begin{equation*}\begin{split}
	&\alpha_i = \frac{(\kappa|\alpha^{\vee}_i)}{n}\quad (i=0,\ldots,2n-1),\quad
	\alpha_{2n} = \rho_1 - \rho_2 + \frac{\kappa_0-\kappa_{2n}}{n},\\
	&\alpha_{2n+1} = -\rho_1 + \rho_2 + \frac{1+\kappa_{2n-1}-\kappa_{2n}}{n},\quad
	\eta = \sum_{j=1}^{n}\frac{\rho_1+\kappa_{2j-2}-\kappa_{2j-1}}{n}.
\end{split}\end{equation*}
Then they satisfy the Painlev\'{e} system $P_{(n+1,n+1)}$.
\end{theorem}

% Section A.2
\subsection{Partition $(2n,1)$}

The Heisenberg subalgebra $\mathfrak{s}_{(2n,1)}$ of $\mathfrak{g}(A_{2n}^{(1)})$ is generated by elements $\{\Lambda_k,\bar{\Lambda}_k\}_{k\in\mathbb{N}}$ with
\begin{equation*}
	\Lambda_1 = \sum_{i=1}^{2n-1}e_i + e_{2n,1},\quad
	\Lambda_2 = \sum_{i=1}^{2n-2}e_{i,1} + e_{2n-1,2} + e_{2n,2}.
\end{equation*}
The grading operator $\vartheta_{(2n,1)}$ is given by
\begin{equation*}
	\vartheta_{(2n,1)} = 2nd + \sum_{i=1}^{2n}\frac{i}{n}\left(2n+1-i-\frac{2}{2n+1}\right)\alpha^{\vee}_i.
\end{equation*}
It implies a $\mathbb{Z}$-gradation of type $(1,\ldots,1,0)$.

The similarity reduction \eqref{Eq:DS_SR_ODE} for $\mathbf{n}=(2n,1)$ is given by
\begin{equation}\label{Eq:SR_2n1}
	[\mathcal{M},\partial_1-B_1] = 0,
\end{equation}
with
\begin{equation*}
	\mathcal{M} = \kappa - \sum_{i=0}^{2n}\varphi_ie_i - \varphi_{0,1}e_{0,1} - \varphi_{2n-1,1}e_{2n-1,1} - \varphi_{2n,1}e_{2n,1} - 2\Lambda_2,
\end{equation*}
where
\begin{equation*}
	\kappa = \vartheta_{(2n,1)} - \sum_{j=0}^{2n}\kappa_j\alpha^{\vee}_j.
\end{equation*}
Note that $\mathcal{M}$ is defined as a $\mathfrak{b}_{+}$-component of $\exp(\mathrm{ad}W)(\vartheta_{(2n,1)}-\rho_1\sum_{i=1}^{2n}i\alpha^{\vee}_i-t_1\Lambda_1-2\Lambda_2)$.

\begin{theorem}
Under the system \eqref{Eq:SR_2n1}, we set
\begin{equation*}\begin{split}
	&t = -\frac{n}{4}t_1^2,\quad
	q_i = \frac{1}{nt_1}\sum_{j=1}^{i}\varphi_{2j},\quad
	p_i = \frac{t_1\varphi_{2i+1}}{4},\quad (i=1,\ldots,n-1),\\
	&q_n = \frac{1}{nt_1}\left(\frac{2\varphi_{2n}}{\varphi_{2n-1,1}}+\sum_{i=1}^{n-1}\varphi_{2i}\right),\quad
	p_n = \frac{t_1\varphi_{0,1}\varphi_{2n-1,1}}{8},
\end{split}\end{equation*}
and
\begin{equation*}\begin{split}
	&\alpha_i = \frac{(\kappa|\alpha^{\vee}_{i+1})}{2n}\quad (i=0,\ldots,2n-1),\\
	&\alpha_{2n} = -\rho_1 - \frac{\kappa_0-\kappa_{2n}}{2n},\quad
	\alpha_{2n+1} = \rho_1 + \frac{1-\kappa_0+\kappa_1}{2n}.
\end{split}\end{equation*}
Then they satisfy the Painlev\'{e} system $P_{(2n+2)}$.
\end{theorem}

% Section A.3
\subsection{Partition $(2n-1,1)$}

The Heisenberg subalgebra $\mathfrak{s}_{(2n-1,1)}$ of $\mathfrak{g}(A_{2n-1}^{(1)})$ is generated by elements $\{\Lambda_k,\bar{\Lambda}_k\}_{k\in\mathbb{N}}$ with
\begin{equation*}
	\Lambda_1 = \sum_{i=1}^{2n-2}e_i + e_{2n-1,1},\quad
	\Lambda_2 = \sum_{i=1}^{2n-3}e_{i,1} + e_{2n-2,2} + e_{2n-1,2}.
\end{equation*}
The grading operator $\vartheta_{(2n-1,1)}$ is given by
\begin{equation*}
	\vartheta_{(2n-1,1)} = (2n-1)d + \sum_{i=1}^{2n-1}\frac{i}{2}\left(2n-i-\frac{1}{n}\right)\alpha^{\vee}_i.
\end{equation*}
It implies a $\mathbb{Z}$-gradation of type $(1,\ldots,1,0)$.

The similarity reduction \eqref{Eq:DS_SR_ODE} for $\mathbf{n}=(2n-1,1)$ is given by
\begin{equation}\label{Eq:SR_2n-11}
	[\mathcal{M},\partial_1-B_1] = 0,
\end{equation}
with
\begin{equation*}
	\mathcal{M} = \kappa - \sum_{i=0}^{2n-1}\varphi_ie_i - \varphi_{0,1}e_{0,1} - \varphi_{2n-2,1}e_{2n-2,1} - \varphi_{2n-1,1}e_{2n-1,1} - 2\Lambda_2,
\end{equation*}
where
\begin{equation*}
	\kappa = \vartheta_{(2n-1,1)} - \sum_{j=0}^{2n-1}\kappa_j\alpha^{\vee}_j.
\end{equation*}
Note that $\mathcal{M}$ is defined as a $\mathfrak{b}_{+}$-component of $\exp(\mathrm{ad}W)(\vartheta_{(2n-1,1)}-\rho_1\sum_{i=1}^{2n-1}i\alpha^{\vee}_i-t_1\Lambda_1-2\Lambda_2)$.

\begin{theorem}
Under the system \eqref{Eq:SR_2n-11} for $\mathbf{n}=(2n-1,1)$, we set
\begin{equation*}\begin{split}
	&t = \frac{\sqrt{2(2n-1)}t_1}{2},\quad
	q_1 = -\frac{\varphi_{0,1}\varphi_{2n-2,1}}{2\sqrt{2(2n-1)}},\quad
	p_1 = \frac{2\varphi_{2n-1}}{\sqrt{2(2n-1)}\varphi_{2n-2,1}},\\
	&q_i = -\frac{\varphi_{0,1}\varphi_{2n-2,1}}{2\sqrt{2(2n-1)}} - \sum_{j=1}^{i-1}\frac{\varphi_{2n-2j}}{\sqrt{2(2n-1)}},\quad
	p_i = \frac{\varphi_{2n-2i+1}}{\sqrt{2(2n-1)}}\quad (i=2,\ldots,n),
\end{split}\end{equation*}
and
\begin{equation*}\begin{split}
	&\alpha_0 = \rho_1 + \frac{1-\kappa_0+\kappa_1}{2n-1},\quad
	\alpha_1 = -\rho_1 - \frac{\kappa_0-\kappa_{2n-1}}{2n-1},\\
	&\alpha_i = \frac{(\kappa|\alpha^{\vee}_{2n-i+1})}{2n-1}\quad (i=2,\ldots,2n-1)
\end{split}\end{equation*}
Then they satisfy the Painlev\'{e} system $P_{(2n+1)}$.
\end{theorem}

% Section B
\section{Symmetric form}

In this appendix, we rewrite the $(2n+2)$-th order Hamiltonian system \eqref{Eq:FS} with \eqref{Eq:FS_Ham} into more simple expression.
Such expression is convenient when we consider the system $P_{(n+1,n+1)}$ from a viewpoint of the affine Weyl group symmetry.

Under the system \eqref{Eq:FS} with \eqref{Eq:FS_Ham}, we set
\begin{equation*}
	x_i = \frac{w_{2i+1}}{t_1^{i-n+\rho_1+\kappa_{2n+1}-\kappa_0}},\quad
	y_i = \frac{t_1^{i-n+\rho_1+\kappa_{2n+1}-\kappa_0}\varphi_{2i+1}}{n+1}\quad (i=0,\ldots,n).
\end{equation*}
Then the Poisson structure is given by
\begin{equation*}
	\{y_i,x_j\} = \delta_{i,j}\quad (i,j=0,\ldots,n).
\end{equation*}
And those variables satisfy the following Hamiltonian system:
\begin{equation}\label{Eq:FS_Sym}
	\frac{dx_i}{dt} = \{H,x_i\},\quad
	\frac{dy_i}{dt} = \{H,y_i\}\quad (i=0,\ldots,n),
\end{equation}
with a Hamiltonian
\begin{equation}\begin{split}\label{Eq:FS_Sym_Ham}
	H &= \frac{1}{t}\sum_{i=0}^{n}\left\{\frac{1}{2}x_i^2y_i^2-\sum_{j=2i+2}^{2n+1}\alpha_jx_iy_i+\sum_{j=0}^{i-1}x_i(x_iy_i+\alpha_{2i+1})y_j\right\}\\
	&\quad + \frac{1}{1-t}\sum_{i=0}^{n}\sum_{j=0}^{n}x_i(x_iy_i+\alpha_{2i+1})y_j,
\end{split}\end{equation}
and a relation
\begin{equation}\label{Eq:FS_Sym_var_rel}
	\sum_{i=0}^{n}x_iy_i + \eta = 0.
\end{equation}
In the following, we call the system \eqref{Eq:FS_Sym} with \eqref{Eq:FS_Sym_Ham} and \eqref{Eq:FS_Sym_var_rel} {\it a symmetric form}.
Recall that
\begin{equation*}
	t = \frac{1}{t_1^{n+1}},\quad
	\alpha_i = \frac{(\kappa|\alpha^{\vee}_i)}{n+1}\quad (i=0,\ldots,2n+1),\quad
	\eta = \sum_{j=0}^{n}\frac{\rho_1+\kappa_{2j}-\kappa_{2j+1}}{n+1}.
\end{equation*}
Note that the canonical coordinates for $P_{(n+1,n+1)}$ are given by
\begin{equation*}
	q_i = \frac{tx_{i-1}}{x_n},\quad
	p_i = \frac{y_{i-1}x_n}{t}\quad (i=1,\ldots,n).
\end{equation*}

The symmetric form is given as the compatibility condition of the Lax pair \eqref{Eq:SR_n+1n+1_Lax} with matrix components
\begin{equation*}\begin{split}
	&\varphi_{2i-1} = y_i\quad (i=1,\ldots,n+1),\\
	&\varphi_{2i} = x_i - x_{i-1}\quad (i=1,\ldots,n),\quad
	\varphi_0 = x_0 - \frac{x_n}{t},
\end{split}\end{equation*}
and
\begin{equation*}\begin{split}
	&u_1 = tv_{2n+1}x_0 - \varepsilon_1,\quad
	u_{2i+1} = v_{2i-1}x_i - \varepsilon_{2i+1}\quad (i=1,\ldots,n),\\
	&u_{2i+2} = v_{2i+1}x_i + \varepsilon_{2n+2}\quad (i=0,\ldots,n),\\
	&v_{2i+1} = \frac{1}{t-1}\left(\sum_{j=0}^{i}y_j+t\sum_{j=i+1}^{n}y_j\right)\quad (i=0,\ldots,n),\\
	&v_{2i} = x_{i-1}\quad (i=1,\ldots,n),\quad
	v_0 = \frac{x_n}{t}.
\end{split}\end{equation*}

The symmetric form is invariant under the actions of the birational canonical transformations $r_0,\ldots,r_{2n+1}$ and $\pi$.
Their actions on the dependent variables are explicitly described as
\[\begin{split}
	&r_0(x_j) = t^{-\alpha_0}x_j,\quad
	r_0(y_j) = t^{\alpha_0}\left(y_j+\frac{\alpha_0}{x_n-tx_0}\{x_n-tx_0,y_j\}\right),\\
	&r_{2i}(x_j) = x_j,\quad
	r_{2i}(y_j) = y_j + \frac{\alpha_{2i}}{x_{i-1}-x_i}\{x_{i-1}-x_i,y_j\}\quad (i=1,\ldots,n),\\
	&r_{2i+1}(x_j) = x_j + \frac{\alpha_{2i+1}}{y_i}\{y_i,x_j\},\quad
	r_{2i+1}(y_j) = y_j\quad (i=0,\ldots,n-1),\\
	&r_{2n+1}(x_j) = t^{\alpha_{2n+1}}\left(x_j+\frac{\alpha_{2n+1}}{y_n}\{y_n,x_j\}\right),\quad
	r_{2n+1}(y_j) = t^{-\alpha_{2n+1}}y_j,
\end{split}\]
for $j=0,\ldots,n$ and
\begin{equation*}\begin{split}
	&\pi(x_i) = \frac{t^{\alpha_0}}{t-1}\left(\sum_{j=0}^{i}y_j+t\sum_{j=i+1}^{n}y_j\right)\quad (i=0,\ldots,n),\\
	&\pi(y_i) = t^{-\alpha_0}(x_i-x_{i+1})\quad (i=0,\ldots,n-1),\quad
	\pi(y_n) = t^{-\alpha_0}(x_n-tx_0).
\end{split}\end{equation*}
We do not describe the actions on the independent variable and the parameters here; see Section \ref{Sec:CP6}.

% Acknowledgement
\section*{Acknowledgement}
The author is grateful to Professors Laszlo Feh\'{e}r, Kenta Fuji, Saburo Kakei, Masatoshi Noumi, Hidetaka Sakai, Teruhisa Tsuda and Yasuhiko Yamada for valuable discussions and advices.

% References

\end{document}